# SCALING LIMITS FOR RANDOM FIELDS WITH LONG-RANGE DEPENDENCE

By Ingemar Kaj,[1] Lasse Leskelä,[2] Ilkka Norros[1]
and Volker Schmidt

*Uppsala University, Helsinki University of Technology, VTT Technical Research Centre and University of Ulm*

This paper studies the limits of a spatial random field generated by uniformly scattered random sets, as the density $\lambda$ of the sets grows to infinity and the mean volume $\rho$ of the sets tends to zero. Assuming that the volume distribution has a regularly varying tail with infinite variance, we show that the centered and renormalized random field can have three different limits, depending on the relative speed at which $\lambda$ and $\rho$ are scaled. If $\lambda$ grows much faster than $\rho$ shrinks, the limit is Gaussian with long-range dependence, while in the opposite case, the limit is independently scattered with infinite second moments. In a special intermediate scaling regime, there exists a nontrivial limiting random field that is not stable.

**1. Introduction.** Fractional Brownian motion often appears as a renormalized limit of independent superpositions of long-memory stochastic processes that are used in physics and other application areas, such as telecommunications and finance. Observing fractional Brownian motion in the limit typically requires rescaling of two model parameters, and switching the order of taking the double limit may lead to approximations with completely different statistical properties [15]. This was also the conclusion of [12], who studied data traffic models with heavy tails, and identified conditions for convergence to fractional Brownian motion and stable Lévy motion in terms

Received June 2005; revised March 2006.

[1]Supported by a basic grant in Applied Mathematics 2003–2005 of the Swedish Foundation for Strategic Research.

[2]Supported by the Mittag–Leffler Institute and the Finnish Graduate School in Stochastics.

of relative scaling speeds of model parameters. This type of results been refined in studies of a special scaling regime that leads to limit processes that are not stable [2, 4, 5, 6].

This paper extends the above trichotomy into a multidimensional context by studying renormalized limits of a spatial random field generated by independently and uniformly scattered random sets in $\mathbb{R}^d$. Viewing the random field as a random linear functional indexed by suitable test functions or test measures, we find different limits for the model as the mean density $\lambda$ of the random sets grows to infinity and the mean volume $\rho$ of the sets tends to zero. If $\lambda$ grows much faster than $\rho$ shrinks, the model with heavy tails converges to a Gaussian self-similar random field with long-range dependence, which in the symmetric case corresponds to fractional Gaussian noise with Hurst parameter $H > 1/2$. In the opposite case where $\rho$ shrinks to zero very rapidly, the limit has infinite second moments and no spatial dependence. We also describe a special intermediate scaling regime that leads to limits that are not stable. In dimension one, these findings correspond to results obtained earlier for a stochastic process known as the infinite source Poisson model or the $M/G/\infty$ model [6, 8, 9, 12].

The outline of the paper is as follows. In Section 2 we construct the model and introduce a functional analytic approach suitable for asymptotic analysis of the random fields. In Section 3 we discuss the different scaling regimes and state the main limit theorems. Section 4 contains a discussion on the statistical properties of the limits, and Section 5 concludes with the proofs.

## 2. Random grain model.
Let $C$ be a bounded measurable set in $\mathbb{R}^d$ such that $|C| = 1$ and $|\partial C| = 0$, where $\partial C$ is the boundary of $C$ and $|\cdot|$ denotes the Lebesgue measure. The building blocks of the model are the sets $x + v^{1/d}C$, called grains, where $x$ is a point in $\mathbb{R}^d$ and $v > 0$ is the volume of the grain. Our goal is to study the mass distribution generated by a family of grains $X_j + (\rho V_j)^{1/d}C$ with random locations $X_j$ and random volumes $\rho V_j$, $j = 1, 2, \ldots$. We assume that $X_j$ are uniformly distributed in the space according to a Poisson random measure with mean density $\lambda > 0$, and that $V_j$ are independent copies of a positive random variable $V$ with $\mathrm{E}V = 1$, also independent of the locations $X_j$. Hence, the scalar $\rho > 0$ equals the mean grain volume. The random field $J_{\lambda,\rho}(x)$ is defined as the number of grains covering $x$,

$$J_{\lambda,\rho}(x) = \#\{j : x \in X_j + (\rho V_j)^{1/d}C\},$$

and we let

(1) $$J_{\lambda,\rho}(A) = \sum_j |A \cap (X_j + (\rho V_j)^{1/d}C)|$$

be the cumulative mass induced by the grains to a measurable set $A$.

2.1. *The random grain field as a random linear functional.* The random variables $J_{\lambda,\rho}(A)$ with $A$ ranging over all measurable sets constitute a set-indexed random function. More generally, we can view $J_{\lambda,\rho}$ as a random functional by replacing the measure $|A \cap \cdot|$ in (1) by an arbitrary positive measure $\phi$,

$$J_{\lambda,\rho}(\phi) = \sum_j \phi(X_j + (\rho V_j)^{1/d} C).$$

Denote by $F(v)$ the probability distribution of $V$. Then the grain volumes $\rho V_j$ are distributed according to $F_\rho(v) = F(v/\rho)$, and $J_{\lambda,\rho}(\phi)$ can be conveniently described as a stochastic integral with respect to a Poisson random measure $N_{\lambda,\rho}(dx, dv)$ on $\mathbb{R}^d \times \mathbb{R}_+$ with intensity measure $\lambda \, dx F_\rho(dv)$,

$$(2) \qquad J_{\lambda,\rho}(\phi) = \int_{\mathbb{R}^d} \int_{\mathbb{R}_+} \phi(x + v^{1/d} C) N_{\lambda,\rho}(dx, dv).$$

To study the linear structure of $J_{\lambda,\rho}$ in a natural way, we do not want to restrict to positive measures. Let $M^1$ be the linear space of signed measures $\phi$ on $\mathbb{R}^d$ with finite total variation $\|\phi\|_1 < \infty$. When $\phi \in M^1$, we see by writing $\phi(A) = \int_A \phi(dx)$ and changing the order of integration that

$$\int_{\mathbb{R}^d} \int_{\mathbb{R}_+} |\phi(x + v^{1/d} C)| \lambda \, dx F_\rho(dv) \le \lambda \rho \|\phi\|_1 < \infty,$$

so the stochastic integral on the right-hand side of (2) converges in probability for all $\phi \in M^1$ [7].

To each function $\phi \in L^1$, one can uniquely associate a signed measure $\tilde{\phi} \in M^1$ defined by $\tilde{\phi}(dx) = \phi(x) \, dx$. We will identify the space $L^1$ with its image in $M^1$ under the map $\phi \mapsto \tilde{\phi}$, so that $L^1 \subset M^1$. Accordingly, when $\phi \in L^1$, we will from now on use the same symbol $\phi$ to signify both the function $\phi(x)$ and the measure $\phi(dx)$. Moreover, if $A$ is a measurable set with $|A| < \infty$, we identify $A$ with the indicator function $1_A \in L^1 \subset M^1$. Note that then $J_{\lambda,\rho}(1_A) = J_{\lambda,\rho}(A)$ agrees with (1).

Denote by $B_r$ the open ball centered at the origin with radius $r$. Then we see that $J_{\lambda,\rho}$ has long-range dependence in the sense that

$$(3) \qquad \lim_{r \to \infty} |\operatorname{Cov}(J_{\lambda,\rho}(B_1), J_{\lambda,\rho}(B_r \setminus B_1))| = \infty$$

if and only if $\mathrm{E} V^2 = \infty$. To verify this, note first that the left-hand side of (3) can be written as

$$(4) \qquad \int_{\mathbb{R}_+} \int_{\mathbb{R}^d} |B_1 \cap (x + v^{1/d} C)| |B_1^c \cap (x + v^{1/d} C)| \lambda \, dx F_\rho(dv)$$

using the covariance formula (10) below. Because $|B_1^c \cap (x + v^{1/d} C)| \le v$ and $\int_{\mathbb{R}^d} |B_1 \cap (x + v^{1/d} C)| \, dx = |B_1| v$, expression (4) is bounded above by

$\lambda\rho^2|B_1|EV^2$. For the other direction, $|B_1^c \cap (x + v^{1/d}C)| \geq v - |B_1|$ implies that $|B_1^c \cap (x + v^{1/d}C)| \geq v/2$ for all $v \geq 2|B_1|$, so that (4) is finite only if $EV^2$ is finite. Observe also that in dimension one, long-range dependence as defined in (3) is equivalent to

$$\sum_{n=-\infty}^{\infty} |\operatorname{Cov}(Y_0, Y_n)| = \infty,$$

where $Y_n = J_{\lambda,\rho}((n, n+1])$ is the discretized version of $J_{\lambda,\rho}$.

2.2. *Riesz energy of signed measures.* To study the limiting behavior of $J_{\lambda,\rho}(\phi)$ as $\lambda \to \infty$ and $\rho \to 0$, we need to impose some more regularity for the measures $\phi \in M^1$. The following subspaces of $M^1$ will turn out to be useful. For $\alpha \in (0, 1)$, let us define

$$M^\alpha = \left\{ \phi \in M^1 : \int_{\mathbb{R}^d} \int_{\mathbb{R}^d} \frac{|\phi|(dx)|\phi|(dy)}{|x-y|^{(1-\alpha)d}} < \infty \right\},$$

where $|\phi|$ is the total variation measure of $\phi$, and for $\phi, \psi \in M^\alpha$, let

$$(5) \qquad \langle \phi, \psi \rangle_\alpha = c_{\alpha,d} \int_{\mathbb{R}^d} \int_{\mathbb{R}^d} \frac{\phi(dx)\psi(dy)}{|x-y|^{(1-\alpha)d}},$$

where

$$(6) \qquad c_{\alpha,d} = \pi^{(\alpha-1/2)d} \Gamma\left(\frac{(1-\alpha)d}{2}\right) \Big/ \Gamma\left(\frac{\alpha d}{2}\right).$$

A classical result in potential theory states that $\langle \phi, \psi \rangle_\alpha$ is an inner product on the vector space $M^\alpha$ [10]. We denote the corresponding norm by $\|\phi\|_\alpha = \langle \phi, \phi \rangle_\alpha^{1/2}$. The quantity $\|\phi\|_\alpha^2$ is often called the Riesz energy of $\phi$. The following proposition describes how the spaces $M^\alpha$ can be ordered.

PROPOSITION 1. *For all $0 < \alpha_1 < \alpha_2 < 1$,*

$$L^1 \cap L^2 \subset M^{\alpha_1} \subset M^{\alpha_2} \subset M^1.$$

REMARK. Let $\mathcal{S}'$ be the space of tempered distributions on $\mathbb{R}^d$, and denote the Fourier transform by $\mathcal{F} : \mathcal{S}' \to \mathcal{S}'$. Then $M^\alpha \subset \mathcal{S}'$ and

$$(7) \qquad \langle \phi, \psi \rangle_\alpha = \int_{\mathbb{R}^d} \mathcal{F}\phi(x) \overline{\mathcal{F}\psi(x)} |x|^{-\alpha d} dx$$

for all $\phi, \psi \in M^\alpha$ ([10], Section VI.1). Equation (7) shows that $\mathcal{F}$ maps $M^\alpha$ isometrically into $L_\alpha^2$, the space of square integrable functions with respect to $|x|^{-\alpha d} dx$. It is also known that $\mathcal{F}(M^\alpha)$ is dense in $L_\alpha^2$, so the Plancherel theorem implies that the closure of $M^\alpha$ with respect to the norm $\|\phi\|_\alpha$ equals $\mathcal{F}^{-1}(L_\alpha^2)$, which is called the space of distributions with finite Riesz energy [10].

2.3. *Integrals with respect to centered Poisson random measures.* Recall that the centered integral $\int f(dN - d\eta)$ of a nonrandom function $f$ with respect to a Poisson random measure $N$ with intensity measure $\eta$ may be defined even for functions that are not $\eta$-integrable. It is known that $\int f(dN - d\eta)$ exists as a limit in probability if and only if

$$(8) \qquad\qquad \int (|f| \wedge f^2) \, d\eta < \infty,$$

in which case the distribution of $\int f(dN - d\eta)$ is characterized by

$$(9) \qquad\qquad \mathrm{E} \exp\Big( i \int f(dN - d\eta) \Big) = \exp \int (\Psi \circ f) \, d\eta,$$

where $\Psi(v) = e^{iv} - 1 - iv$ for $v \in \mathbb{R}$ [7]. Moreover,

$$(10) \qquad\qquad \mathrm{E}\Big( \int f(dN - d\eta) \Big)\Big( \int g(dN - d\eta) \Big) = \int fg \, d\eta,$$

when $f$ and $g$ are square integrable with respect to $\eta$.

## 3. Scaling behavior and main results.

3.1. *Scaling behavior of the random grain model.* We will next study the limiting behavior of $J_{\lambda,\rho}(\phi)$ as the mean grain density $\lambda$ grows to infinity and the mean grain volume $\rho$ shrinks to zero. When the grain volume distribution has finite variance, the following central limit theorem shows that the centered and renormalized version of $J_{\lambda,\rho}$ converges to white Gaussian noise.

THEOREM 1. *Let $C$ be a bounded set with $|C| = 1$ and $|\partial C| = 0$, and assume $\mathrm{E}V^2 < \infty$. Then as $\lambda \to \infty$ and $\rho \to 0$, the following limit holds in the sense of finite-dimensional distributions of random functionals indexed by $L^1 \cap L^2$:*

$$\frac{J_{\lambda,\rho}(\phi) - \mathrm{E}J_{\lambda,\rho}(\phi)}{\rho(\lambda \mathrm{E}V^2)^{1/2}} \longrightarrow W(\phi),$$

*where $W$ is the centered Gaussian random linear functional on $L^2$ with*

$$(11) \qquad\qquad \mathrm{E}W(\phi)W(\psi) = \int_{\mathbb{R}^d} \phi(x)\psi(x) \, dx.$$

However, our main focus will be on the model where the volume distribution is heavy-tailed with infinite variance. Hence, we will from now assume that the distribution $F(v)$ of the normalized volume $V$ has a regularly varying tail of index $\gamma \in (1, 2)$, that is,

$$(12) \qquad\qquad \lim_{v \to \infty} \frac{\bar{F}(av)}{\bar{F}(v)} = a^{-\gamma} \qquad \text{for all } a > 0,$$

where $\bar{F}(v) = 1 - F(v)$. This implies that $\mathrm{E}V^2 = \infty$. Let us denote $f(\rho) \sim g(\rho)$, if $f(\rho)/g(\rho) \to 1$ as $\rho \to 0$. Then (12) implies that the scaled volume distribution $F_\rho(v) = F(v/\rho)$ satisfies

$$\bar{F}_\rho(v) \sim \bar{F}_\rho(1)v^{-\gamma} \qquad \text{as } \rho \to 0,$$

and by Karamata's theorem [see formula (22) in Section 5.1], the expected number of grains with volume larger than one that cover the origin equals

$$\iint_{\{(x,v)\,:\,0\in x+v^{1/d}C,v>1\}} \lambda\,dx\,F_\rho(dv) = \lambda\int_1^\infty vF_\rho(dv) \sim \frac{\lambda\bar{F}_\rho(1)}{1-\gamma^{-1}}.$$

Consequently, we distinguish the following three scaling regimes:

$$\text{large-grain scaling } \lambda\bar{F}_\rho(1) \to \infty,$$
$$\text{intermediate scaling } \lambda\bar{F}_\rho(1) \to \sigma_0 > 0,$$
$$\text{small-grain scaling } \lambda\bar{F}_\rho(1) \to 0.$$

The regular variation of $\bar{F}(v)$ implies that the relations $\lambda \sim (1/\rho)^{\gamma+\varepsilon}$ and $\lambda \sim (1/\rho)^{\gamma-\varepsilon}$ for some $\varepsilon > 0$ belong to large-grain and small-grain regimes, respectively, while in the critical intermediate scaling regime, the size of $\lambda$ is roughly proportional to $(1/\rho)^\gamma$.

Under large-grain scaling, the number of grains that are big enough to carry statistical dependence over macroscopic distances grows to infinity. Hence, the limit of $J_{\lambda,\rho}$ in this case is expected to have long-range spatial dependence. In the opposite case of small-grain scaling, no grains survive that are big enough to represent substantial dependence over spatial distances, so the small-grain limit of $J_{\lambda,\rho}$ should have very weak dependence over space. The intermediate scaling regime is a blend of the two other, with a balanced mix of large grains providing long-range spatial dependence and small grains generating nontrivial random variations on short distances.

3.2. *Main results.* The following theorem justifies the heuristics in Section 3.1. Let $\gamma \in (1,2)$, and recall that the independently scattered $\gamma$-stable random measure with unit skewness and Lebesgue control measure is the random linear functional $\Lambda_\gamma(\phi) = \int \phi(x)\Lambda_\gamma(dx)$ on $L^\gamma$ characterized by

$$(13) \qquad \mathrm{E}e^{i\Lambda_\gamma(\phi)} = \exp\left(-\sigma_\phi^\gamma\left(1 - i\beta_\phi\tan\left(\frac{\pi\gamma}{2}\right)\right)\right),$$

where $\sigma_\phi = \|\phi\|_\gamma$ and $\beta_\phi = \|\phi\|_\gamma^{-\gamma}(\|\phi_+\|_\gamma^\gamma - \|\phi_-\|_\gamma^\gamma)$, and where $\phi_+ = \max(\phi,0)$ and $\phi_- = -\min(\phi,0)$. For an alternative equivalent definition of $\Lambda_\gamma$ as a set-indexed random function, see [14].

THEOREM 2. *Let $C$ be a bounded set with $|C| = 1$ and $|\partial C| = 0$, and assume that $V$ has a regularly varying tail with exponent $\gamma \in (1, 2)$. Let $\alpha \in (0, 2 - \gamma)$. Then the following three limits hold in the sense of finite-dimensional distributions of random functionals as $\lambda \to \infty$ and $\rho \to 0$:*

(i) *(Large-grain scaling) If $\lambda \bar{F}_\rho(1) \to \infty$, then*

$$\frac{J_{\lambda,\rho}(\phi) - \mathrm{E} J_{\lambda,\rho}(\phi)}{(\gamma \lambda \bar{F}_\rho(1))^{1/2}} \longrightarrow W_{\gamma,C}(\phi), \qquad \phi \in M^\alpha,$$

*where $W_{\gamma,C}$ is the centered Gaussian random linear functional on $M^{2-\gamma}$ with $\mathrm{E} W_{\gamma,C}(\phi) W_{\gamma,C}(\psi) = \int \int \phi(dx) K_{\gamma,C}(x - y) \psi(dy)$ and*

$$K_{\gamma,C}(x) = \int_0^\infty |(v^{-1/d}x + C) \cap C| v^{-\gamma} \, dv. \tag{14}$$

(ii) *(Intermediate scaling) If $\lambda \bar{F}_\rho(1) \to \sigma_0 > 0$, then*

$$J_{\lambda,\rho}(\phi) - \mathrm{E} J_{\lambda,\rho}(\phi) \longrightarrow J^*_{\gamma,C}(\phi_\sigma), \qquad \phi \in M^\alpha,$$

*where $J^*_{\gamma,C}$ is defined on $M^{2-\gamma}$ as a centered integral with respect to the Poisson random measure $N_\gamma(dx, dv)$ on $\mathbb{R}^d \times \mathbb{R}_+$ with intensity measure $dx \, v^{-\gamma-1} \, dv$,*

$$J^*_{\gamma,C}(\phi) = \int_{\mathbb{R}^d} \int_0^\infty \phi(x + v^{1/d}C)(N_\gamma(dx, dv) - dx \, v^{-\gamma-1} \, dv), \tag{15}$$

*and where $\phi_\sigma$ is defined by $\phi_\sigma(A) = \phi(\sigma A)$ with $\sigma = (\gamma \sigma_0)^{1/((\gamma-1)d)}$.*

(iii) *(Small-grain scaling) If $\lambda \bar{F}_\rho(1) \to 0$, then*

$$\frac{J_{\lambda,\rho}(\phi) - \mathrm{E} J_{\lambda,\rho}(\phi)}{c_\gamma (1/\bar{F}_\rho)^\leftarrow(\gamma\lambda)} \longrightarrow \Lambda_\gamma(\phi), \qquad \phi \in L^1 \cap L^2,$$

*where $\Lambda_\gamma$ is the independently scattered $\gamma$-stable random measure on $\mathbb{R}^d$ with Lebesgue control measure and unit skewness, $(1/\bar{F}_\rho)^\leftarrow(u) = \inf\{v : 1/\bar{F}_\rho(v) \geq u\}$ is the quantile function of $F_\rho$, and*

$$c_\gamma = \left( -\frac{\Gamma(2 - \gamma)}{\gamma(\gamma - 1)} \cos\left(\frac{\pi\gamma}{2}\right) \right)^{-1/\gamma}. \tag{16}$$

3.3. *Role of symmetry and randomly oriented grains.* Fix a parameter $H \in (1/2, 1)$, and let $W_H$ be the centered Gaussian random linear functional on $M^{2H-1}$ with

$$\mathrm{E} W_H(\phi) W_H(\psi) = c_{2H-1,d} \int_{\mathbb{R}^d} \int_{\mathbb{R}^d} \frac{\phi(dx)\psi(dy)}{|x - y|^{(2-2H)d}} = \langle \phi, \psi \rangle_{2H-1}, \tag{17}$$

*where $\langle \phi, \psi \rangle_{2H-1}$ is the Riesz inner product defined by (5) and $c_{2H-1,d}$ is given by (6). When $C$ is symmetric around the origin so that $\theta C = C$ for* all rotations $\theta$ of $\mathbb{R}^d$, the rotation invariance of the Lebesgue measure shows that the covariance kernel of the large-grain limit $W_{\gamma,C}$ satisfies

$$K_{\gamma,C}(x) = K_{\gamma,C}(|x|e_1) = K_{\gamma,C}(e_1)|x|^{-(\gamma-1)d},$$

where $e_1$ is an arbitrary fixed unit vector in $\mathbb{R}^d$. For symmetric grains $C$, it hence follows that $W_{\gamma,C}$ equals $cW_H$ in the sense of finite-dimensional distributions, where $H = (3-\gamma)/2$ and $c = c_{2H-1,d}^{-1}K_{\gamma,C}(e_1)^{1/2}$.

We can also study the limit behavior for a slightly modified model where the grains have independent and uniform random orientations. To define this model, let $d\theta$ be the Haar measure on the compact group $SO(d)$ of rotations in $\mathbb{R}^d$, and let $N_{\lambda,\rho}(dx, dv, d\theta)$ be a Poisson random measure on $\mathbb{R}^d \times \mathbb{R}_+ \times SO(d)$ with intensity measure $\lambda\,dx F_\rho(dv)\,d\theta$. Then

$$\widetilde{J}_{\lambda,\rho}(\phi) = \int_{\mathbb{R}^d} \int_{\mathbb{R}_+} \int_{SO(d)} \phi(x + v^{1/d}\theta C) N_{\lambda,\rho}(dx, dv, d\theta)$$

defines the analogue of $J_{\lambda,\rho}$ with randomly rotated grains $\theta C$ [compare with definition (2) in Section 2.1]. Because $d\theta$ is a probability measure on the compact group $SO(d)$ that is not scaled during $\lambda \to \infty$ and $\rho \to 0$, the following result can be verified by copying the proof of Theorem 2. Note that the shape of $C$ reduces into a constant $c$ in the large-grain limit below.

THEOREM 3. *Under the assumptions of Theorem 2, the following three limits hold in the sense of finite-dimensional distributions of random functionals as $\lambda \to \infty$ and $\rho \to 0$:*

(i) *(Large-grain scaling) If $\lambda \bar{F}_\rho(1) \to \infty$, then*

$$\frac{\widetilde{J}_{\lambda,\rho}(\phi) - \mathrm{E}\widetilde{J}_{\lambda,\rho}(\phi)}{(\gamma\lambda\bar{F}_\rho(1))^{1/2}} \longrightarrow cW_H(\phi), \qquad \phi \in M^\alpha,$$

*where $W_H$ is the Gaussian random linear functional defined in (17) with $H = (3-\gamma)/2$ and*

$$c = c_{2H-1,d}^{-1} \left( \int_{SO(d)} \int_0^\infty |(v^{-1/d}\theta e_1 + C) \cap C|v^{-\gamma}\,dv\,d\theta \right)^{1/2}.$$

(ii) *(Intermediate scaling) If $\lambda \bar{F}_\rho(1) \to \sigma_0 > 0$, then*

$$\widetilde{J}_{\lambda,\rho}(\phi) - \mathrm{E}\widetilde{J}_{\lambda,\rho}(\phi)$$
$$\longrightarrow \int_{\mathbb{R}^d} \int_0^\infty \int_{SO(d)} \phi_\sigma(x + v^{1/d}\theta C)(N_\gamma(dx, dv, d\theta) - dx\,v^{-\gamma-1}\,dv\,d\theta),$$

*$\phi \in M^\alpha$, where $N_\gamma(dx, dv, d\theta)$ is a Poisson random measure on $\mathbb{R}^d \times \mathbb{R}_+ \times SO(d)$ with intensity $dx\,v^{-\gamma-1}dv\,d\theta$ and $\phi_\sigma$ is as in Theorem 2.*

(iii) (Small-grain scaling) If $\lambda \bar{F}_\rho(1) \to 0$, then

$$\frac{\widetilde{J}_{\lambda,\rho}(\phi) - \mathrm{E}\widetilde{J}_{\lambda,\rho}(\phi)}{c_\gamma(1/\bar{F}_\rho)^\leftarrow(\gamma\lambda)} \longrightarrow \Lambda_\gamma(\phi), \qquad \phi \in L^1 \cap L^2,$$

where $\Lambda_\gamma$ and $c_\gamma$ are as in Theorem 2.

## 4. Statistical properties of the limits.

4.1. *Properties of the large-grain limit.* A change of variables shows that the covariance kernel of $W_{\gamma,C}$ given by (14) scales according to $K_{\gamma,C}(ax) = a^{-(\gamma-1)d}K_{\gamma,C}(x)$ for $a > 0$. Hence, for $H = (3-\gamma)/2$,

$$\mathrm{E}W_{\gamma,C}(\phi_s)W_{\gamma,C}(\psi_s) = s^{2(1-H)d}\mathrm{E}W_{\gamma,C}(\phi)W_{\gamma,C}(\psi),$$

where the dilatation $\phi_s$ of the measure $\phi$ is defined for $s > 0$ by

$$(18) \qquad\qquad \phi_s(A) = \phi(sA).$$

Thus, $W_{\gamma,C}$ is self-similar in the sense that $W_{\gamma,C}(\phi_s)$ and $s^{(1-H)d}W_{\gamma,C}(\phi)$ have the same finite-dimensional distributions on $M^{2H-1}$ for all $s > 0$.

To study the autocovariance properties of $W_{\gamma,C}$ over long spatial ranges, note first that

$$\lim_{r\to\infty} |\mathrm{Cov}(W_{\gamma,C}(B_1), W_{\gamma,C}(B_r \setminus B_1))| = \int_{B_1}\int_{B_1^c} K_{\gamma,C}(x-y)\,dy\,dx.$$

By changing the order of integration,

$$\int_{B_1^c} K_{\gamma,C}(x-y)\,dy = \int_0^\infty \int_{x-v^{1/d}C} |B_1^c \cap (y+v^{1/d}C)|\,dy v^{-\gamma-1}\,dv,$$

and because $|B_1^c \cap (y+v^{1/d}C)| \ge v - |B_1|$ for all $v$, we see that the inner integral on the right-hand side above is greater than or equal to $v^2/2$ for all $v \ge 2|B_1|$. Hence, for all $x$,

$$\int_{B_1^c} K_{\gamma,C}(x-y)\,dy \ge \frac{1}{2}\int_{2|B_1|}^\infty v^{-\gamma+1}\,dv,$$

which is infinite for $\gamma \in (1,2)$. From this, we conclude that

$$\lim_{r\to\infty} |\mathrm{Cov}(W_{\gamma,C}(B_1), W_{\gamma,C}(B_r \setminus B_1))| = \infty,$$

which means that $W_{\gamma,C}$ has long-range dependence in the sense of (3).

The symmetric large-grain limit $W_H$ can be represented in terms of the white Gaussian noise $W$ defined in Theorem 1, when $W$ is viewed as an independently scattered Gaussian random measure with Lebesgue control measure as in [14].

PROPOSITION 2. *For $H \in (1/2, 1)$, the random linear functional $W_H$ on $M^{2H-1}$ equals*

$$(19) \qquad W_H(\phi) = c_{H-1/2,d} \int_{\mathbb{R}^d} \int_{\mathbb{R}^d} \frac{\phi(dy)}{|x-y|^{(3/2-H)d}} W(dx),$$

*in the sense of finite-dimensional distributions, where $c_{H-1/2,d}$ is given by* (6).

Specializing to dimension one, we see that for $\phi, \psi \in L^1 \cap M^{2H-1}$,

$$\mathrm{E} W_H(\phi) W_H(\psi) = c_{2H-1,1} \int_{-\infty}^{\infty} \int_{-\infty}^{\infty} \phi(s) \psi(t) |t-s|^{2H-2} \, ds \, dt,$$

from which we recognize that $W_H(\phi)$ equals a constant multiple of the stochastic integral of $\phi$ with respect to fractional Brownian motion with Hurst parameter $H > 1/2$ [3]. The functional $W_H$ may thus be viewed as a natural extension of fractional Gaussian noise [11] into multidimensional parameter spaces. Moreover, choosing $\phi = 1_{[0,t]}$ in (19) yields the well-balanced representation of fractional Brownian motion ([14], Section 7.2.1).

4.2. *Properties of the intermediate limit.* Using (10) and changing the order of integration,

$$\mathrm{E} J_{\gamma,C}^*(\phi) J_{\gamma,C}^*(\psi) = \int_{\mathbb{R}^d} \int_{\mathbb{R}_+} \phi(x + v^{1/d}C) \psi(x + v^{1/d}C) \, dx \, v^{-\gamma-1} \, dv$$

$$= \int_{\mathbb{R}^d} \int_{\mathbb{R}^d} \phi(dx) K_{\gamma,C}(x-y) \phi(dy),$$

which shows that $J_{\gamma,C}^*$ and $W_{\gamma,C}$ share the same second order statistical structure. Especially, this implies that $J_{\gamma,C}^*$ has long-range dependence in the sense of (3).

We will next show that $J_{\gamma,C}^*$ is not self-similar by assuming the contrary and deriving a contradiction. Assume that $J_{\gamma,C}^*(\phi_s) = a_s J_{\gamma,C}^*(\phi)$ in distribution for all $s > 0$, where $\phi_s(A) = \phi(sA)$ as before. Then the self-similarity of $W_{\gamma,C}$ implies that $a_s = s^{(\gamma-1)d/2}$, because $\mathrm{E} J_{\gamma,C}^*(\phi_s)^2 = \mathrm{E} W_{\gamma,C}(\phi_s)^2$. A change of variables shows that

$$\int_{\mathbb{R}^d} \int_{\mathbb{R}_+} \Psi(\phi_s(x + v^{1/d}C)) \, dx \, v^{-\gamma-1} \, dv$$

$$= s^{(\gamma-1)d} \int_{\mathbb{R}^d} \int_{\mathbb{R}_+} \Psi(\phi(x + v^{1/d}C)) \, dx \, v^{-\gamma-1} \, dv.$$

Comparing this with the characteristic functional of $J_{\gamma,C}^*$ given by (9) and denoting $t = s^{(\gamma-1)d/2}$ allows us to conclude that

$$\int_{\mathbb{R}^d} \int_{\mathbb{R}_+} \frac{\Psi_{re}(t\phi(x + v^{1/d}C))}{t^2} \, dx \, v^{-\gamma-1} \, dv$$

(20)
$$= \int_{\mathbb{R}^d} \int_{\mathbb{R}_+} \Psi_{\mathrm{re}}(\phi(x + v^{1/d}C)) \, dx \, v^{-\gamma-1} \, dv,$$

where $\Psi_{\mathrm{re}}$ denotes the real part of $\Psi$. Because $|\Psi_{\mathrm{re}}(v)| \le 2$ for all $v$, it follows that the integrand on the left-hand side of (20) converges to zero as $t \to \infty$. Moreover, $|\Psi_{re}(v)| \le v^2/2$ (Lemma 1) implies that this sequence of integrands is bounded from above by the function $\phi(x + v^{1/d}C)^2/2$, which is integrable with respect to $v^{-\gamma-1} \, dv \, dx$, as verified in (36) below. Hence, by dominated convergence, the left-hand side of (20) converges to zero as $t \to \infty$. When $\phi$ is chosen so that the right-hand side of (20) is nonzero, this is a contradiction.

A similar reasoning can be used to verify that $J_{\gamma,C}^*$ is not stable. However, the sum of $n$ independent copies of $J_{\gamma,C}^*$ has the same finite-dimensional distributions as $\phi \mapsto J_{\gamma,C}^*(\phi_s)$, where $s = n^{1/((\gamma-1)d)}$. This property is called aggregate-similarity in [4].

### 4.3. *Properties of the small-grain limit.*

Note that when $\phi$ is a function in $L^1 \cap L^2 \subset M^1$, the dilatation $\phi_s$ defined in (18) becomes
$$\phi_s(A) = \int_{sA} \phi(x) \, dx = \int_A s^d \phi(sx) \, dx,$$

so for functions, $\phi_s(x) = s^d \phi(sx)$. Inspection of the characteristic functional (13) of $\Lambda_\gamma$ shows that $\Lambda_\gamma(\phi_s)$ and $s^{(1-1/\gamma)d}\Lambda_\gamma(\phi)$ have the same finite-dimensional distributions, so $\Lambda_\gamma$ is self-similar. The random functional $\Lambda_\gamma$ can also be represented by
$$\Lambda_\gamma(\phi) = c_\gamma^{-1} \int_{\mathbb{R}^d} \int_0^\infty v\phi(x)(N_\gamma(dx, dv) - dx \, v^{-\gamma-1} \, dv),$$

where $N_\gamma(dx, dv)$ is the Poisson random measure appearing in (15) and $c_\gamma$ is given by (16); see [14]. Specializing to dimension one, we remark that the stochastic process
$$t \mapsto \Lambda_\gamma(1_{[0,t]}) = \int_0^t \Lambda_\gamma(ds)$$

is the centered $\gamma$-stable Lévy motion with unit skewness.

### 5. Proofs.

By definition, $J_{\lambda,\rho}$ together with the four limit fields defined in Theorem 1 and Theorem 2 are linear in the sense that for all test measures $\phi_1, \dots, \phi_n$ and all scalars $a_1, \dots, a_n$,
$$J_{\lambda,\rho}(a_1\phi_1 + \cdots + a_n\phi_n) = a_1 J_{\lambda,\rho}(\phi_1) + \cdots + a_n J_{\lambda,\rho}(\phi_n)$$

almost surely. Hence, convergence of the finite-dimensional distributions of the centered and renormalized version of $J_{\lambda,\rho}$ is equivalent to the convergence of the one-dimensional distributions. Recall from (9) that for $b > 0$, the characteristic functional of $(J_{\lambda,\rho}(\phi) - \mathrm{E}J_{\lambda,\rho}(\phi))/b$ is given by

$$
(21) \quad
\begin{aligned}
&\mathrm{E}\exp\left(i\frac{J_{\lambda,\rho}(\phi) - \mathrm{E}J_{\lambda,\rho}(\phi)}{b}\right) \\
&\quad = \exp\int_{\mathbb{R}^d}\int_{\mathbb{R}_+}\Psi\left(\frac{\phi(x + v^{1/d}C)}{b}\right)\lambda F_\rho(dv)\,dx,
\end{aligned}
$$

where $\Psi(v) = e^{iv} - 1 - iv$. The following lemma summarizes the properties of $\Psi$ that are needed in proving the theorems of the paper.

LEMMA 1.  *The function* $\Psi(v) = e^{iv} - 1 - iv$ *satisfies*

$$|\Psi(v) - \Psi(u)| \le (2|v - u| \wedge |v^2 - u^2|/2)$$

*for all* $u, v \in \mathbb{R}$. *Moreover, for all* $v \in \mathbb{R}$,

$$|\Psi(v)| \le (2|v| \wedge v^2/2) \quad and \quad |\Psi(v) + v^2/2| \le (v^2 \wedge |v|^3/6).$$

PROOF.  Observe first that $\Psi(v + 2n\pi) - \Psi(u + 2n\pi) = \Psi(v) - \Psi(u)$ for all integers $n$. Hence, in proving the first inequality, we may without loss of generality assume that $u$ and $v$ are nonnegative. Moreover, by symmetry, it is enough to consider the case $u \le v$. For $0 \le u \le v$, we have

$$|\Psi(v) - \Psi(u)| \le \int_u^v |e^{is} - 1|\,ds \le \int_u^v (2 \wedge s)\,ds,$$

because $|e^{is} - 1| \le (2 \wedge |s|)$ for all $s$. This proves the first inequality. The second inequality follows from the first by setting $u = 0$. Further, the third inequality follows from the second because $\Psi(v) + v^2/2 = i\int_0^v \Psi(s)\,ds$.  □

Before going to the proofs of the main theorems, we first introduce some preliminary results on regular variation (Section 5.1) and on maximal functions (Section 5.2). Proposition 1 is proved in Section 5.3, while in Section 5.4 we develop the key results on the regularity of the characteristic functional (21). Sections 5.5–5.8 contain the proofs of the main theorems, and Section 5.9 concludes with the proof of Proposition 2.

5.1. *Regular variation.*  Let $F$ be a probability distribution on $\mathbb{R}_+$ with a regularly varying tail of exponent $\gamma > 0$, and let $0 < p < \gamma < q$. Then using integration by parts and Karamata's theorem ([1], Theorem 1.5.11) it follows that as $a \to \infty$,

$$
(22) \quad \int_{(a,\infty)} v^p F(dv) \sim \frac{\gamma}{\gamma - p}\bar{F}(a)a^p,
$$

$$
(23) \quad \int_{[0,a]} v^q F(dv) \sim \frac{\gamma}{q - \gamma}\bar{F}(a)a^q,
$$

where $\bar{F}(a) = 1 - F(a)$. The next two lemmas summarize the theory on regular variation that is later used to analyze the distribution of the normalized volume $V$.

LEMMA 2. *Let $F$ be a probability distribution on $\mathbb{R}_+$ with a regularly varying tail of exponent $\gamma > 0$, and define the scaled distribution $F_\rho$ for $\rho > 0$ by $F_\rho(v) = F(v/\rho)$. Assume that $f(v)$ is a continuous function on $\mathbb{R}_+$ such that, for some $0 < p < \gamma < q$,*

$$\limsup_{v \to \infty} v^{-p}|f(v)| < \infty \quad and \quad \limsup_{v \to 0} v^{-q}|f(v)| < \infty.$$

*Then*

$$\int_{\mathbb{R}_+} f(v) F_\rho(dv) \sim \bar{F}_\rho(1) \int_0^\infty f(v) \gamma v^{-\gamma-1} \, dv \qquad as \; \rho \to 0.$$

PROOF. Fix a constant $a \in (0,1)$. Then (22) implies that, for all $v_0 > 0$,

$$\int_{(v_0,\infty)} v^p F_\rho(dv) \sim \bar{F}_\rho(1)\gamma \int_{v_0}^\infty v^{-\gamma-1+p} \, dv,$$

which shows that the finite measures $\bar{F}_\rho(1)^{-1} v^p F_\rho(dv)$ restricted to $(a,\infty)$ converge weakly to the finite measure $\gamma v^{-\gamma-1+p} \, dv$ on $(a,\infty)$ as $\rho \to 0$. Because the function $v^{-p} f(v)$ is continuous and bounded on $(a,\infty)$, this implies

$$(24) \qquad \bar{F}_\rho(1)^{-1} \int_{(a,\infty)} f(v) F_\rho(dv) \to \int_{(a,\infty)} f(v)\gamma v^{-\gamma-1} \, dv.$$

Moreover, the second assumption on $f$ implies that $|f(v)| \le cv^q$ for all $v \in [0,1]$, so (23) implies that

$$\left| \bar{F}_\rho(1)^{-1} \int_{[0,a]} f(v) F_\rho(dv) - \int_0^a f(v)\gamma v^{-\gamma-1} \, dv \right|$$

$$\le c\bar{F}_\rho(1)^{-1} \int_{[0,a]} v^q F_\rho(dv) + c \int_0^a v^q \gamma v^{-\gamma-1} \, dv \sim 2c\frac{\gamma}{q-\gamma} a^{q-\gamma}.$$

The claim now follows because the right-hand side above can be made arbitrarily small by choosing $a$ small enough, and because (24) holds for all $a \in (0,1)$. □

LEMMA 3. *Let $F$ and $F_\rho$ be defined as in Lemma 2, and let $f_\rho(v)$ be a family of measurable functions on $\mathbb{R}_+$. Assume that for some $0 < p < \gamma < q$, either*

$$(25) \qquad \begin{aligned} &\limsup_{\rho \to 0} \sup_{v > a} v^{-p} \bar{F}_\rho(1)|f_\rho(v)| = 0 && for \; all \; a > 0, \\ &\bar{F}_\rho(1)|f_\rho(v)| \le cv^q && for \; all \; \rho, v, \end{aligned}$$

*or*

$$\lim_{\rho \to 0} \sup_{a \le v} v^{-q} \bar{F}_\rho(1)|f_\rho(v)| = 0 \qquad \text{for all } a > 0,$$

(26)

$$\bar{F}_\rho(1)|f_\rho(v)| \le cv^p \qquad \text{for all } \rho, v.$$

*Then*

$$\lim_{\rho \to 0} \int_{\mathbb{R}_+} f_\rho(v) F_\rho(dv) = 0.$$

PROOF. Assume that the functions $f_\rho(v)$ satisfy the conditions (25) and fix $a > 0$. Denote $c_a(\rho) = \sup_{v > a} v^{-p} \bar{F}_\rho(1)|f_\rho(v)|$. Then by (22),

$$\int_{(a,\infty)} |f_\rho(v)| F_\rho(dv) \le c_a(\rho) \bar{F}_\rho(1)^{-1} \int_{(a,\infty)} v^p F_\rho(dv) \sim c_a(\rho) \frac{\gamma}{\gamma - p} a^{p-\gamma},$$

and by (23),

$$\int_{(0,a]} |f_\rho(v)| F_\rho(dv) \le c \bar{F}_\rho(1)^{-1} \int_{[0,a]} v^q F_\rho(dv) \sim c \frac{\gamma}{q - \gamma} a^{q-\gamma}.$$

Because $c_a(\rho) \to 0$ as $\rho \to 0$, the two above bounds imply that

$$\limsup_{\rho \to 0} \int_{\mathbb{R}_+} |f_\rho(v)| F_\rho(dv) \le c \frac{\gamma}{q - \gamma} a^{q-\gamma}.$$

Since this is true for all $a > 0$, the claim follows by letting $a \to 0$. The proof under assumption (26) is analogous. □

5.2. *Maximal functions.* Let $C$ be a bounded measurable set in $\mathbb{R}^d$ with $|C| = 1$. If $\phi$ is a locally integrable function, define the averages $m_\phi(x, v)$ by

(27)
$$m_\phi(x, v) = v^{-1} \int_{x + v^{1/d}C} \phi(y) \, dy,$$

and let $\phi_*$ be the maximal function of $\phi$ given by

(28)
$$\phi_*(x) = \sup_{v > 0} v^{-1} \int_{x + v^{1/d}C} |\phi(y)| \, dy.$$

The following lemma summarizes the known facts about maximal functions that are used to find integrable upper bounds for the characteristic functional of $J_{\lambda, \rho}$ in the proofs of Theorem 1 and Theorem 2(iii).

LEMMA 4. *Let $C$ be a bounded measurable set in $\mathbb{R}^d$ with $|C| = 1$:*

(i) *If $\phi \in L^1$, then $\lim_{v \to 0} m_\phi(x, v) = \phi(x)$ for almost all $x$.*

(ii) *If $\phi \in L^1$, then $\phi_*(x) < \infty$ for almost all $x$.*

(iii) *If $\phi \in L^p$ for some $p > 1$, then $\phi_* \in L^p$.*

Proof. Let $B$ be the open ball centered at the origin with unit volume, and fix $a > 0$ such that $C \subset a^{1/d} B$. Assume first that $\phi \in L^1$. Then

$$v^{-1} \int_{x+v^{1/d} C} |\phi(y) - \phi(x)| \, dy \le a(av)^{-1} \int_{x+(av)^{1/d} B} |\phi(y) - \phi(x)| \, dy$$

for all $x \in \mathbb{R}^d$ and $v > 0$. The right-hand side tends to zero as $v \to 0$, because almost all points $x \in \mathbb{R}^d$ are Lebesgue points of $\phi$ ([13], Theorem 7.7). Thus, the first claim is valid.

Next, let

$$\phi^*(x) = \sup_{v > 0} v^{-1} \int_{x+v^{1/d} B} |\phi(y)| \, dy$$

be the Hardy–Littlewood maximal function of $\phi$. Then $C \subset a^{1/d} B$ implies that $\phi_*(x) \le a\phi^*(x)$ for all $x$. The second claim now follows, because for $\phi \in L^1$, $\phi^*(x) < \infty$ almost everywhere ([13], Theorem 7.4). The third claim follows directly from the Hardy–Littlewood maximal theorem ([13], Theorem 8.18), which states that if $\phi \in L^p$ for some $p > 1$, then $\phi^* \in L^p$. □

5.3. *Proof of Proposition* 1. Assume $0 < \alpha_1 < \alpha_2 < 1$. When $\phi \in L^1 \cap L^2$, then denoting $D = \{(x, y) : |x - y| \le 1\}$, we see that

$$(29) \qquad \int\int_{D^c} \frac{|\phi(x)||\phi(y)|}{|x - y|^{(1-\alpha_1)d}} \, dx \, dy \le \|\phi\|_1^2.$$

Moreover, writing

$$\int\int_D \frac{|\phi(x)||\phi(y)|}{|x - y|^{(1-\alpha_1)d}} \, dx \, dy = \int\int_D \frac{|\phi(x)|}{|x - y|^{(1-\alpha_1)d/2}} \frac{|\phi(y)|}{|x - y|^{(1-\alpha_1)d/2}} \, dx \, dy$$

and applying Hölder's inequality, we see that

$$\int\int_D \frac{|\phi(x)||\phi(y)|}{|x - y|^{(1-\alpha_1)d}} \, dx \, dy \le \int\int_D \frac{\phi(x)^2}{|x - y|^{(1-\alpha_1)d}} \, dx \, dy$$

$$= \|\phi\|_2^2 \int_{\{x \,:\, |x| \le 1\}} \frac{dx}{|x|^{(1-\alpha_1)d}}.$$

This bound together with (29) shows that $\phi \in M^{\alpha_1}$, so the first inclusion has been shown.

To verify the second inclusion, note that $|x - y|^{(1-\alpha_2)d} \ge |x - y|^{(1-\alpha_1)d}$ on $D$, and $|x - y|^{(1-\alpha_2)d} \ge 1$ on $D^c$. Thus,

$$\int_{\mathbb{R}^d} \int_{\mathbb{R}^d} \frac{|\phi|(dx)|\phi|(dy)}{|x - y|^{(1-\alpha_2)d}} \le \int\int_D \frac{|\phi|(dx)|\phi|(dy)}{|x - y|^{(1-\alpha_1)d}} \, dx \, dy + \int\int_{D^c} |\phi|(dx)||\phi|(dy),$$

which is finite for $\phi \in M^{\alpha_1}$. □

5.4. *Regularity properties of the characteristic functional.* In this section we prove the key continuity and boundedness properties of the characteristic functional of $J_{\lambda,\rho}$ that are required for the asymptotical analysis of the model. We start with a continuity property of the Lebesgue measure. Denote the symmetric difference of sets $A$ and $B$ by

$$A \triangle B = (A \setminus B) \cup (B \setminus A).$$

LEMMA 5. *Let $C$ be a bounded measurable set in $\mathbb{R}^d$ such that $|\partial C| = 0$. Then*

$$\lim_{r \to 1} |C \triangle rC| = 0.$$

PROOF. If $y \neq 0$ belongs to the interior of $C$, then $1_{rC}(y) = 1_C(y/r)$ converges to 1 as $r \to 1$. Moreover, because $|\partial C| = 0$, it follows that $\lim_{r \to 1} 1_{rC}(y) = 1$ for almost all $y \in C$. Hence, the dominated convergence theorem implies

$$\lim_{r \to 1} |C \cap rC| = \lim_{r \to 1} \int_C 1_{rC}(y)\, dy = |C|,$$

and by writing

$$|C \triangle rC| = |C| - |C \cap rC| + |rC| - |C \cap rC|$$
$$= (1 + r^d)|C| - 2|C \cap rC|,$$

we see that the claim is valid. □

LEMMA 6. *Let $C$ be a bounded measurable set in $\mathbb{R}^d$ such that $|\partial C| = 0$. Then for each $\phi \in M^1$, the functions*

$$v \mapsto \int_{\mathbb{R}^d} \phi(x + v^{1/d}C)^2\, dx \quad and \quad v \mapsto \int_{\mathbb{R}^d} \Psi(\phi(x + v^{1/d}C))\, dx$$

*are continuous on $\mathbb{R}_+$.*

PROOF. Define for $u, v \geq 0$,

$$d(u,v) = \int_{\mathbb{R}^d} |\phi(x + v^{1/d}C) - \phi(x + u^{1/d}C)|\, dx.$$

Then, using $|v^2 - u^2| = |u + v||u - v|$, we see that

$$(30) \qquad \int_{\mathbb{R}^d} |\phi(x + v^{1/d}C)^2 - \phi(x + u^{1/d}C)^2|\, dx \leq 2\|\phi\|_1 d(u,v),$$

while $|\Psi(v) - \Psi(u)| \leq 2|v - u|$ (Lemma 1) implies

$$(31) \qquad \int_{\mathbb{R}^d} |\Psi(\phi(x + v^{1/d}C)) - \Psi(\phi(x + u^{1/d}C))|\, dx \leq 2d(u,v).$$

Next, observe that

$$|\phi(x + v^{1/d}C) - \phi(x + u^{1/d}C)| \le |\phi|((x + u^{1/d}C)\triangle(x + v^{1/d}C))$$
$$= |\phi|(x + (u^{1/d}C\triangle v^{1/d}C)).$$

Because $\int |\phi|(x + A)\,dx = \|\phi\|_1|A|$ for all measurable sets $A$, the above inequality implies that

$$d(u, v) \le \|\phi\|_1|u^{1/d}C\triangle v^{1/d}C|.$$

Moreover, because $|u^{1/d}C\triangle v^{1/d}C| = u|C\triangle(v/u)^{1/d}C|$, it follows using Lemma 5 that $\lim_{u\to v} d(u, v) = 0$ for all $v \ge 0$. This fact together with the bounds (30) and (31) completes the proof. $\square$

LEMMA 7. *Let $C$ be a bounded measurable set in $\mathbb{R}^d$, and assume $\phi \in M^\alpha$ with $\alpha \in (0, 1)$. Then there exists a constant $c$ such that, for all $v \ge 0$,*

$$\int_{\mathbb{R}^d} \phi(x + v^{1/d}C)^2\,dx \le c(v \wedge v^{2-\alpha}).$$

PROOF. Let $\phi$ be a measure in $M^\alpha$ with $\alpha \in (0, 1)$. Because the total variation measure $|\phi|$ also belongs to $M^\alpha$, we may assume without loss of generality that $\phi$ is positive. Then by changing the order of integration, we see that

$$(32) \qquad \int_{\mathbb{R}^d} \phi(x + v^{1/d}C)^2\,dx \le \|\phi\|_1 \int_{\mathbb{R}^d} \phi(x + v^{1/d}C)\,dx = |C|\|\phi\|_1^2 v.$$

Next, let $a$ be large enough such that $C \subset aB$, where $B$ is the open unit ball. Then

$$\phi(x + v^{1/d}C)^2 \le \phi(x + (av)^{1/d}B)^2,$$

so again by changing the order of integration, we see that

$$(33) \qquad \begin{aligned} &\int_{\mathbb{R}^d} \phi(x + v^{1/d}C)^2\,dx \\ &\le \int_{\mathbb{R}^d}\int_{\mathbb{R}^d} |(x - (av)^{1/d}B)\cap(y - (av)^{1/d}B)|\phi(dx)\phi(dy). \end{aligned}$$

Let $e_1$ be an arbitrary fixed unit vector in $\mathbb{R}^d$. Because the Lebesgue measure is rotation and translation invariant, we see that, for all $x$ and $y$,

$$\begin{aligned} |(x - (av)^{1/d}B)\cap(y - (av)^{1/d}B)| &= |(|x - y|e_1 + (av)^{1/d}B)\cap(av)^{1/d}B| \\ &\le av1(|x - y| < 2(av)^{1/d}) \\ &\le \frac{2^{(1-\alpha)d}(av)^{2-\alpha}}{|x - y|^{(1-\alpha)d}}, \end{aligned}$$

where the first inequality holds because the set $(|x - y|e_1 + (av)^{1/d}B) \cap (av)^{1/d}B$ is empty for $(av)^{1/d} \leq |x - y|/2$. Combining the above bound with (33), we get

$$(34) \qquad \int_{\mathbb{R}^d} \phi(x + v^{1/d}C)^2 \, dx \leq 2^{(1-\alpha)d} a^{2-\alpha} \|\phi\|_\alpha^2 v^{2-\alpha}.$$

Combining inequalities (32) and (34) together, we conclude that the claim holds by taking $c = \max(|C| \|\phi\|_1^2, 2^{(1-\alpha)d} a^{2-\alpha} \|\phi\|_\alpha^2)$.   $\square$

5.5. *Proof of Theorem* 1. Let $\phi \in L^1 \cap L^2$ and define $b = \rho(\lambda \mathrm{E}V^2)^{1/2}$. Without loss of generality, choose $\rho$ as the basic model parameter and consider $\lambda$ and $b$ as functions of $\rho$. Because $\phi \in L^1 \cap L^2$,

$$\phi(x + v^{1/d}C) = \int_{x + v^{1/d}C} \phi(y) \, dy = v m_\phi(x, v),$$

where $m_\phi(x, v)$ is the average of $\phi$ defined in (27). Using (21) and the definition of $b$, we thus see that

$$\mathrm{E} \exp\left(i \frac{J_{\lambda,\rho}(\phi) - \mathrm{E}J_{\lambda,\rho}(\phi)}{b}\right) = \exp \int_{\mathbb{R}^d} \int_{\mathbb{R}_+} \lambda \Psi\left(\frac{v m_\phi(x, \rho v)}{(\lambda \mathrm{E}V^2)^{1/2}}\right) F(dv) \, dx.$$

By Lemma 4, $\lim_{\rho \to 0} m_\phi(x, \rho v) = \phi(x)$ for all $v$ and almost all $x$. By Lemma 1, $\Psi(v) = -\frac{1}{2}v^2 + \varepsilon(v)$, where $|\varepsilon(v)| \leq |v|^3/6$, so that

$$(35) \qquad \lim_{\rho \to 0} \lambda \Psi\left(\frac{v m_\phi(x, \rho v)}{(\lambda \mathrm{E}V^2)^{1/2}}\right) = -\frac{v^2 \phi(x)^2}{2\mathrm{E}V^2}.$$

Moreover, letting $\phi_*$ be the maximal function of $\phi$ defined in Lemma 4, the bound $|\Psi(v)| \leq v^2/2$ (Lemma 1) implies that

$$\lambda \Psi\left(\frac{v m_\phi(x, \rho v)}{(\lambda \mathrm{E}V^2)^{1/2}}\right) \leq \frac{v^2 \phi_*(x)^2}{2\mathrm{E}V^2}.$$

Because by Lemma 4 the right-hand side is integrable with respect to $dx F(dv)$, the dominated convergence theorem combined with (35) shows that

$$\lim_{\rho \to 0} \mathrm{E} \exp\left(i \frac{J_{\lambda,\rho}(\phi) - \mathrm{E}J_{\lambda,\rho}(\phi)}{b}\right) = \exp\left(-\frac{1}{2} \int_{\mathbb{R}^d} \phi(x)^2 \, dx\right),$$

where the right-hand side is the characteristic functional of the white Gaussian noise $W$ on $L^2$. Hence, the proof of Theorem 1 is complete.   $\square$

5.6. *Proof of Theorem* 2, *large-grain scaling.* Fix $\gamma \in (1, 2)$, let $\phi \in M^{2-\gamma}$, and assume first that $\phi$ is positive. Then by changing the order of integration, it follows that

$$\int_{\mathbb{R}^d} \phi(x + v^{1/d}C)^2 \, dx = \int_{\mathbb{R}^d} \int_{\mathbb{R}^d} |(x - v^{1/d}C) \cap (y - v^{1/d}C)| \phi(dx) \phi(dy),$$

so by the translation invariance of the Lebesgue measure,

$$
\begin{aligned}
(36) \qquad \int_{\mathbb{R}^d} \int_{\mathbb{R}_+} & \phi(x + v^{1/d}C)^2 v^{-\gamma-1} \, dv \, dx \\
&= \int_{\mathbb{R}^d} \int_{\mathbb{R}^d} \phi(dx) K_{\gamma,C}(x - y)\phi(dy),
\end{aligned}
$$

where $K_{\gamma,C}$ is the covariance kernel defined in (14). Choose $a > 0$ large enough so that $C \subset aB$, where $B$ is the open unit ball. Then letting $e_1$ be an arbitrary unit vector in $\mathbb{R}^d$, we see that

$$
K_{\gamma,C}(x) \le K_{\gamma,aB}(x) = K_{\gamma,aB}(e_1)|x|^{(1-\gamma)d},
$$

so the right-hand side of (36) is bounded by $K_{\gamma,aB}(e_1)\|\phi\|_{2-\gamma}^2$ and hence, finite. Thus, by Fubini's theorem, equation (36) holds also for nonpositive $\phi \in M^{2-\gamma}$. Because the left-hand side of (36) is nonnegative, $\int \int \phi(dx)K_{\gamma,C}(x - y)\psi(dy)$ is a positive definite bilinear form in $M^{2-\gamma} \times M^{2-\gamma}$, and hence, defines the distribution of a centered Gaussian random linear functional on $M^{2-\gamma}$, which we call $W_{\gamma,C}$.

Assume next that $\phi \in M^\alpha \subset M^{2-\gamma}$ for some $\alpha \in (0, 2 - \gamma)$ and let $b = (\gamma\lambda\bar{F}_\rho(1))^{1/2}$. As before, we choose $\rho$ as the basic model parameter and consider $\lambda$ and $b$ as functions of $\rho$. Define the functions $f_\rho$ and $f$ on $\mathbb{R}_+$ by

$$
f_\rho(v) = \int_{\mathbb{R}^d} \Psi\left(\frac{\phi(x + v^{1/d}C)}{b}\right) dx, \qquad f(v) = -\frac{1}{2}\int_{\mathbb{R}^d} \phi(x + v^{1/d}C)^2 \, dx.
$$

By Lemma 6, the function $f(v)$ is continuous, and $|f(v)| \le c(v \wedge v^{2-\alpha})/2$ by Lemma 7. Application of Lemma 2 with $p = 1$ and $q = 2 - \alpha$ hence yields

$$
\int_{\mathbb{R}_+} f(v)F_\rho(dv) \sim \bar{F}_\rho(1) \int_0^\infty f(v)\gamma v^{-\gamma-1} \, dv,
$$

so that by the definition of $b$,

$$
(37) \qquad \lim_{\rho \to 0} \int_{\mathbb{R}_+} f(v)\lambda b^{-2}F_\rho(dv) = \int_0^\infty f(v)\gamma v^{-\gamma-1} \, dv.
$$

Next, define $g_\rho(v) = \lambda f_\rho(v) - \lambda b^{-2}f(v)$, and observe that

$$
g_\rho(v) = \lambda \int_{\mathbb{R}^d} \left( \Psi\left(\frac{\phi(x + v^{1/d}C)}{b}\right) + \frac{1}{2}\left(\frac{\phi(x + v^{1/d}C)}{b}\right)^2 \right) dx.
$$

Because $|\Psi(v) - (-v^2/2)| \le |v|^3/6$ (Lemma 1) and

$$
\int_{\mathbb{R}^d} |\phi(x + v^{1/d}C)|^3 \, dx \le \|\phi\|_1^2 \int_{\mathbb{R}^d} |\phi|(x + v^{1/d}C) \, dx = \|\phi\|_1^3 v,
$$

we see that $|g_\rho(v)| \le \lambda b^{-3}\|\phi\|_1^3 v/6$. Using the definition of $b$, we thus see that

$$
(38) \qquad \bar{F}_\rho(1)v^{-1}|g_\rho(v)| \le \frac{\|\phi\|_1^3}{6\gamma b}
$$

for all $v \geq 0$. Moreover, using $|\Psi(v)| \leq v^2/2$ and Lemma 7, it follows that $|g_\rho(v)| \leq c\lambda b^{-2} v^{2-\alpha}$. Hence,

$$\bar{F}_\rho(1) v^{-(2-\alpha)} |g_\rho(v)| \leq c/\gamma$$

for all $v \geq 0$. The large-grain assumption $\lambda \bar{F}_\rho(1) \to \infty$ implies that $b \to \infty$, so that the right-hand side of (38) tends to zero as $\rho \to 0$. Thus, using Lemma 3 with $p = 1$ and $q = 2 - \alpha$, we conclude that

$$(39) \qquad \lim_{\rho \to 0} \int_{\mathbb{R}_+} g_\rho(v) F_\rho(dv) = 0.$$

Combining (37) and (39), we get

$$\lim_{\rho \to 0} \int_{\mathbb{R}_+} f_\rho(v) \lambda F_\rho(dv) = \int_0^\infty f(v) v^{-\gamma-1} \, dv.$$

In light of (21) and (36), this is equivalent to

$$\lim_{\rho \to 0} \mathrm{E} \exp\left( i \frac{J_{\lambda,\rho}(\phi) - \mathrm{E} J_{\lambda,\rho}(\phi)}{b} \right) = \exp\left( -\frac{1}{2} \int_{\mathbb{R}^d} \int_{\mathbb{R}^d} \phi(dx) K_{\gamma,C}(x) \phi(dy) \right),$$

which completes the proof of Theorem 2(i).  $\square$

5.7. *Proof of Theorem 2, intermediate scaling.* Let $\phi$ be a measure in $M^{2-\gamma}$. In the beginning of the proof of Theorem 2, part (i), we saw that

$$\int_{\mathbb{R}_+} \int_{\mathbb{R}^d} \phi(x + v^{1/d} C)^2 \, dx \, v^{-\gamma-1} \, dv < \infty.$$

Thus, by (8), the stochastic integral

$$\int_{\mathbb{R}^d} \int_{\mathbb{R}_+} \phi(x + v^{1/d} C)(N_\gamma(dx, dv) - dx \, v^{-\gamma-1} \, dv)$$

converges in probability, so the right-hand side of (15) is well-defined on $M^{2-\gamma}$.

Assume next that $\phi \in M^\alpha$ for some $\alpha \in (0, 2 - \gamma)$, and choose again $\rho$ as the basic model parameter and consider $\lambda$ as a function of $\rho$. Define

$$f(v) = \int_{\mathbb{R}^d} \Psi(\phi(x + v^{1/d} C)) \, dx.$$

Note that $f(v)$ is continuous by Lemma 6. Moreover, $|\Psi(v)| \leq v^2/2$ (Lemma 1) together with Lemma 7 shows that $|f(v)| \leq c(v \wedge v^{2-\alpha})/2$. Lemma 2 with $p = 1$ and $q = 2 - \alpha$ thus shows that

$$\int_{\mathbb{R}_+} f(v) \lambda F_\rho(dv) \sim \lambda \bar{F}_\rho(1) \int_0^\infty f(v) \gamma v^{-\gamma-1} \, dv.$$

The characteristic functional formula (21) together with the intermediate scaling assumption $\lambda \bar{F}_\rho(1) \to \sigma_0$ now implies that

$$\lim_{\rho \to 0} \mathrm{E} \exp(i(J_{\lambda,\rho}(\phi) - \mathrm{E}J_{\lambda,\rho}(\phi))) = \exp\left(\gamma \sigma_0 \int_0^\infty f(v) v^{-\gamma-1} \, dv\right).$$

Denoting $\sigma = (\gamma \sigma_0)^{1/((\gamma-1)d)}$ and defining $\phi_\sigma(A) = \phi(\sigma A)$, a change of variables shows that the right-hand side above equals

$$\exp\left(\int_{\mathbb{R}^d} \int_0^\infty \Psi(\phi_\sigma(x + v^{1/d}C)) \, dx \, v^{-\gamma-1} \, dv\right).$$

By (9), this agrees with the characteristic functional of (15), so the proof of Theorem 2(ii) is complete. $\square$

5.8. *Proof of Theorem* 2, *small-grain scaling.* Define $b = (1/\bar{F}_\rho)^{\leftarrow}(\gamma\lambda)$ (for notational convenience, we do not include the constant $c_\gamma$ into $b$). As before, $\lambda$ and $b$ are considered to be functions of $\rho$. Because $\lambda \to \infty$, it follows from $(1/\bar{F}_\rho)^{\leftarrow}(\gamma\lambda) = \rho(1/\bar{F})^{\leftarrow}(\gamma\lambda)$ that $b/\rho \to \infty$, and by Theorem 1.5.12 in [1],

$$(40) \qquad \lambda \bar{F}(b/\rho) \sim \gamma^{-1} \qquad \text{as } \rho \to 0.$$

The small-grain scaling assumption $\lambda \bar{F}_\rho(1) \to 0$ implies that $\lambda \bar{F}_\rho(\varepsilon) \to 0$ for all $\varepsilon > 0$. Observe that for all $\rho$ such that $b \geq \varepsilon$, we have

$$\gamma^{-1} \sim \lambda \bar{F}_\rho(b) \leq \lambda \bar{F}_\rho(\varepsilon).$$

Because the right-hand side above converges to zero, $b$ must eventually become less than $\varepsilon$ as $\rho \to 0$. In other words, $b \to 0$.

Let $\phi \in L^1 \cap L^2$. We start by showing that, for almost all $x$,

$$(41) \qquad \lim_{\rho \to 0} \int_{\mathbb{R}_+} \Psi(v m_\phi(x, bv)) \lambda F_{\rho/b}(dv) = \int_0^\infty \Psi(v\phi(x)) v^{-\gamma-1} \, dv.$$

Observe first that because $\Psi(v)$ is continuous and $|\Psi(v)| \leq (2|v| \wedge v^2/2)$ by Lemma 1, we can apply Lemma 2 with $p = 1$ and $q = 2$ (and with $\rho/b \to 0$ in place of $\rho$) to the function $v \mapsto \Psi(v\phi(x))$ and conclude that

$$\int_{\mathbb{R}_+} \Psi(v\phi(x)) \lambda F_{\rho/b}(dv) \sim \lambda \bar{F}_{\rho/b}(1) \int_0^\infty \Psi(v\phi(x)) \gamma v^{-\gamma-1} \, dv.$$

Using (40), this shows that, for all $x$,

$$(42) \qquad \lim_{\rho \to 0} \int_{\mathbb{R}_+} \Psi(v\phi(x)) \lambda F_{\rho/b}(dv) = \int_0^\infty \Psi(v\phi(x)) v^{-\gamma-1} \, dv.$$

Next, let

$$g_\rho(v) = \Psi(v m_\phi(x, bv)) - \Psi(v\phi(x)).$$

By Lemma 1, $|\Psi(v) - \Psi(u)| \le |v^2 - u^2|/2$. Hence, by using $(v - u)^2 = (v - u)(v + u)$, we see that, for all $a > 0$,

$$\sup_{0 < v \le a} v^{-2} \bar{F}_{\rho/b}(1)|g_\rho(v)| \le (\phi_*(x) + |\phi(x)|) \sup_{0 < v \le a} |m_\phi(x, bv) - \phi(x)|/2.$$

By Lemma 4, the right-hand side tends to zero for almost all $x$, as $\rho \to 0$ (recall that $b$ tends to zero together with $\rho$). Moreover, $|\Psi(v)| \le 2|v|$ shows that

$$v^{-1} \bar{F}_{\rho/b}(1)|g_\rho(v)| \le 2(\phi_*(x) + |\phi(x)|)$$

for all $\rho$ and $v$. Thus, we can now use Lemma 3 with $p = 1$ and $q = 2$ (and $\rho/b \to 0$ in place of $\rho$) to conclude that

$$\lim_{\rho \to 0} \int_{\mathbb{R}_+} g_\rho(v) F_{\rho/b}(dv) = 0.$$

Combining this with (42) now shows the validity of (41).

Property (41) implies that

$$
\begin{aligned}
(43) \quad & \lim_{\rho \to 0} \int_{\mathbb{R}^d} \int_{\mathbb{R}_+} \Psi(vm_\phi(x, bv)) \lambda F_{\rho/b}(dv)\, dx \\
& = \int_{\mathbb{R}^d} \int_0^\infty \Psi(v\phi(x)) v^{-\gamma - 1}\, dv\, dx,
\end{aligned}
$$

provided we can take the limit in (43) inside the $dx$-integral. To justify this interchange of the limit and the integral, choose a small enough $\varepsilon > 0$ such that $\gamma \in (1 + \varepsilon, 2 - \varepsilon)$. Then $|\Psi(v)| \le 2\min(|v|, v^2) \le 2\min(|v|^{\gamma - \varepsilon}, |v|^{\gamma + \varepsilon})$ and $|m_\phi(x, bv)| \le \phi_*(x)$ imply that

$$|\Psi(vm_\phi(x, bv))| \le 2(\phi_*(x)^{\gamma - \varepsilon} + \phi_*(x)^{\gamma + \varepsilon})(v^{\gamma - \varepsilon} \wedge v^{\gamma + \varepsilon}).$$

Moreover, by Lemma 2 (with $\rho/b$ in place of $\rho$),

$$\int_{\mathbb{R}_+} (v^{\gamma - \varepsilon} \wedge v^{\gamma + \varepsilon}) \lambda F_{\rho/b}(dv) \sim \lambda \bar{F}_{\rho/b}(1) \int_0^\infty (v^{\gamma - \varepsilon} \wedge v^{\gamma + \varepsilon}) \gamma v^{-\gamma - 1}\, dv,$$

so by (40), we see that the integral on the left-hand side converges to $2\varepsilon^{-1}$, and hence, becomes eventually less than $1 + 2\varepsilon^{-1}$ as $\rho \to 0$. Thus, for all $\rho$ small enough,

$$\int_{\mathbb{R}_+} |\Psi(vm_\phi(x, bv))| \lambda F_{\rho/b}(dv) \le 2(1 + 2\varepsilon^{-1})(\phi_*(x)^{\gamma - \varepsilon} + \phi_*(x)^{\gamma + \varepsilon}).$$

By Lemma 4, the right-hand side above is $dx$-integrable. Thus, the dominated convergence theorem shows the validity of (43). Further, using (21),

$$\lim_{\rho \to 0} \mathrm{E} \exp\left( i \frac{J_{\lambda, \rho}(\phi) - \mathrm{E} J_{\lambda, \rho}(\phi)}{b} \right) = \exp\left( \int_{\mathbb{R}^d} \int_0^\infty \Psi(v\phi(x)) v^{-\gamma - 1}\, dv\, dx \right).$$

By splitting the integration over $\mathbb{R}^d$ into $\{x : \phi(x) \geq 0\}$ and $\{x : \phi(x) < 0\}$ and performing a change of variables, one can verify that the right-hand side above equals

$$\exp(d_\gamma \|\phi_+\|_\gamma^\gamma + \bar{d}_\gamma \|\phi_-\|_\gamma^\gamma),$$

where $\bar{d}_\gamma$ is the complex conjugate of $d_\gamma = \int_0^\infty \Psi(v) v^{-\gamma-1}\, dv$. Moreover,

$$d_\gamma = \frac{\Gamma(2-\gamma)}{\gamma(\gamma-1)} \cos\left(\frac{\pi\gamma}{2}\right)\left(1 - i\tan\left(\frac{\pi\gamma}{2}\right)\right);$$

see Exercise 3.24 in [14]. Comparing the definition of $c_\gamma$ given in (16) with the characteristic functional of $\Lambda_\gamma$ in (13), we conclude that

$$\lim_{\rho \to 0} \mathrm{E} \exp\left(i \frac{J_{\lambda,\rho}(\phi) - \mathrm{E} J_{\lambda,\rho}(\phi)}{b}\right) = \mathrm{E} e^{ic_\gamma \Lambda_\gamma(\phi)},$$

which completes the proof of Theorem 2(iii). $\quad\square$

5.9. *Proof of Proposition* 2. Assume that $\phi$ is a positive measure in $M^{(2H-1)}$, let $\alpha = 2H - 1$, and define

$$f(x) = c_{\alpha/2,d} \int_{\mathbb{R}^d} \frac{\phi(dy)}{|x-y|^{(1-\alpha/2)d}}.$$

Then the composition rule for Riesz kernels ([10] Section 1.1) shows that

$$c_{\alpha/2,d}^2 \int_{\mathbb{R}^d} \frac{1}{|x-y|^{(1-\alpha/2)d}} \frac{1}{|x-y'|^{(1-\alpha/2)d}}\, dx = c_{\alpha,d} \frac{1}{|y-y'|^{(1-\alpha)d}}.$$

Hence, by changing the order of integration, we see that

$$\int_{\mathbb{R}^d} f(x)^2\, dx = c_{\alpha,d} \int_{\mathbb{R}^d} \int_{\mathbb{R}^d} \frac{\phi(dy)\phi(dy')}{|y-y'|^{(1-\alpha)d}} = \langle \phi, \phi \rangle_\alpha,$$

where $\langle \phi, \psi \rangle_\alpha$ is the Riesz inner product defined by (5). Comparing this with (11) shows that the Gaussian random variables on the left and the right-hand side of (19) have the same variance, and thus equal in distribution. Equality of the finite-dimensional distributions follows by linearity. $\quad\square$

I. KAJ
DEPARTMENT OF MATHEMATICS
UPPSALA UNIVERSITY
P.O. BOX 480
S-751 06 UPPSALA
SWEDEN
E-MAIL: ikaj@math.uu.se

I. NORROS
VTT TECHNICAL RESEARCH CENTRE
P.O. BOX 1202
FI-02044 VTT
FINLAND
E-MAIL: ilkka.norros@vtt.fi

L. LESKELÄ
INSTITUTE OF MATHEMATICS
HELSINKI UNIVERSITY OF TECHNOLOGY
P.O. BOX 1100
FI-02015 TKK
FINLAND
E-MAIL: lasse.leskela@iki.fi

V. SCHMIDT
DEPARTMENT OF STOCHASTICS
UNIVERSITY OF ULM
D-89069 ULM
GERMANY
E-MAIL: volker.schmidt@mathematik.uni-ulm.de